\documentclass[12pt]{amsart}
\usepackage{amscd,amssymb,float}
\newtheorem{theorem}{Theorem}[section]
\newtheorem{lemma}[theorem]{Lemma}
\newtheorem{claim}[theorem]{Claim}
\newtheorem{proposition}[theorem]{Proposition}

\theoremstyle{definition}

\theoremstyle{remark}

\theoremstyle{conj}

\numberwithin{equation}{section}

\begin{document}
\title{Convergence of K\"ahler-Ricci flow with integral curvature bound}

\author[F. Fang]{Fuquan Fang}
\thanks{The first author was supported by a
NSF Grant of China and the Capital Normal University}
\address{Department of Mathematics, Capital Normal University,
Beijing, P.R.China}
  \email{ffang@nankai.edu.cn}
\author[Y. Zhang]{Yuguang Zhang}
\address{Department of Mathematics, Capital Normal University,
Beijing, P.R.China  }
 \email{zhangyuguang76@yahoo.com}

\begin{abstract}Let $g(t)$, $t\in [0, +\infty)$, be a solution of the
normalized K\"ahler-Ricci flow  on a compact K\"ahler $n$-manifold
$M$ with $c_{1}(M)>0$ and initial metric  $g (0)\in 2\pi
c_{1}(M)$.
 If there is a constant $C$ independent of $t$ such that
$$
\int_{M}|Rm(g(t))|^{n}dv_{t}\leq C,$$  then, for any
$t_{k}\rightarrow \infty$, a subsequence of $(M, g(t_{k}))$
converges to a compact  orbifold $(X, h)$ with only finite many
singular points $\{q_{j}\}$ in the Gromov-Hausdorff sense, where
$h$ is a  K\"ahler metric on $X\backslash \{q_{j}\}$ satisfying
 the K\"ahler-Ricci soliton equation, i.e.  there is a  smooth function
$f$ such that
$$Ric(h)-h=\nabla\overline{\nabla}f,  \ \ \ \ {\rm and}\it \ \ \ \ \ \ \nabla \nabla f=\overline{\nabla} \overline{\nabla} f=0.
$$
\end{abstract}
\maketitle

\section{Introduction}
On a compact K\"ahler $n$-manifold $M$ with $c_{1}(M)>0$, the
normalized K\"ahler-Ricci flow equation  is
\begin{equation}\label{1.1}\partial_{t}g(t)=-Ric(g(t))+g(t)=\sqrt{-1}\partial \overline{\partial}u_{t},
\end{equation} for a family of K\"ahler metrics $g(t)\in 2\pi c_{1}(M)$, where we identify K\"ahler
 metrics with the K\"ahler forms.   In
\cite{Cao}, it is proved that  a solution $g(t)$ of (\ref{1.1})
exists
 for all times $t\in [0, \infty)$.  Perelman (cf. \cite{ST})) has proved some
 important properties for the solution $g(t)$, $t\in [0, \infty)$,  of
 (\ref{1.1}): there exist constants $C>0$ $\kappa >0$ independent of $t$ such
 that
 \begin{itemize}\label{00}
  \item[(1)] $|R(g(t))|< C$, and  $\text{diam}_{g(t)}(M)<C$,
 \item[(2)] $|u_{t}|_{C^{1}(g(t))} < C$,
  \item[(3)] $(M, g(t))$ is $\kappa$-noncollapsed, i.e. for any
$r< 1$,  if $|R(g(t))|\leq r^{-2}$ on a metric ball
$B_{g(t)}(x,r)$, then
\begin{equation}\label{2.5}\text{Vol}_{g(t)}(B_{g(t)}(x,r))\geq \kappa r^{2n}.
\end{equation}

 \end{itemize} By assuming that  the Ricci curvature is uniformly
 bounded along the flow, Sesum and Tian have proved that, for any
  sequence of times  $t_{k}\longrightarrow \infty$, a subsequence of $(M,
 g(t_{k}+t))$ converges to $(X, g_{\infty}(t))$, where $X$ is
 smooth outside a singular set, and $g_{\infty}(t)$ satisfies the K\"ahler-Ricci soliton
 equation (cf. \cite{ST} and \cite{S2}). In a recent preprint \cite{S3},
 Sesum has proved that $X$ is actually a K\"ahler manifold if $n\geq
 3$, and $g(t)$ satisfies an additional  integral bound  of curvature
 operators.
   The purpose  of this note is to  study  the
 convergence of $(M,
 g(t_{k}))$  by assuming an integral bound  of curvature operators
 instead of the  uniform bound for Ricci curvatures.

\vskip 5mm

\begin{theorem} Let $g(t)$, $t\in [0, +\infty)$, be a solution of the
normalized K\"ahler-Ricci flow (\ref{1.1}) on a compact K\"ahler
$n$-manifold $M$ with $c_{1}(M)>0$ and initial metric  $g (0)\in
2\pi c_{1}(M)$.
 If there is a constant $C$ independent of $t$ such that
$$
\int_{M}|Rm(g(t))|^{n}dv_{t}\leq C,$$  then, for any
$t_{k}\rightarrow \infty$, a subsequence of $(M, g(t_{k}))$
converges to a compact  orbifold $(X, h)$ with only finite many
singular points $\{q_{j}\}$ in the Gromov-Hausdorff sense.
Furthermore, on $X\backslash \{q_{j}\}$,  $h$ is a  K\"ahler
metric satisfying
 the K\"ahler-Ricci soliton equation, i.e.  there is a  smooth function
$f$ such that
$$Ric(h)-h=\nabla\overline{\nabla}f,  \ \ \ \ {\rm and}\it \ \ \ \ \ \ \nabla \nabla f=\overline{\nabla} \overline{\nabla} f=0.
$$
\end{theorem}
\vskip 5mm

Here, we call a topological space $X$ an orbifold if $X$ is a
smooth manifold outside a finite set of singular points
$\{q_{j}\}$, and there is a neighborhood around every singular
point homeomorphic to a cone on a spherical space form
$C(S^{2n-1}/\Gamma)$, $\Gamma\subset SO(2n)$. A metric $h$ on $X$
is a Riemannian metric on $X\backslash \{q_{j}\}$, and, in a local
uniformization $B^{2n}\backslash \{0\}$, $h$ extends to a
$C^{0}$-metric on the ball $B^{2n}$. Note that this definition is
 different from the one in \cite{An2}, which allows several
 spherical cones joint at a single vertex.

 In   \cite{P2}, Perelman claimed that, if $M$ admits a
 K\"ahler-Einstein metric  with positive scalar curvature, then a
 solution of (\ref{1.1}) converges to  the K\"ahler-Einstein metric. In
 the case of $M$ admitting a shrinking  K\"ahler-Ricci soliton, Tian and Zhu
 obtained the same result in \cite{TZ}. However, there are K\"ahler
 manifolds with  $c_{1}>0$  admitting no K\"ahler-Einstein
 metrics and no holomorphic vector field, thus no K\"ahler-Ricci
 soliton (\cite{T2}). In this situation, solutions  of (\ref{1.1})
 will develop singularities when times tend to infinity by
 Theorem 1.1 and \cite{S2}.

If $n=2$, the $L^{2}$-norm of the curvature operator of a K\"ahler
metric is uniformly bounded by terms of the first and the  second
Chern class and its K\"ahler class (c.f. \cite{S} and \cite{CCT}).
In this case, Sesum has claimed a strong version of Theorem 1.1
basing on an unpublished work on the  K\"ahler-Ricci flow due to
Tian (c.f. \cite{S}). \vskip 5mm

\begin{theorem}[Sesum, Tian] Let $g(t)$, $t\in [0, +\infty)$, be a solution of the
normalized K\"ahler-Ricci flow (\ref{1.1}) on a compact K\"ahler
 surface  $M$ with $c_{1}(M)>0$ and initial metric  $g (0)\in 2\pi
c_{1}(M)$. Then, for any $t_{k}\rightarrow \infty$, a subsequence
of $(M, g(t_{k}))$ converges to an orbifold $(X, h)$ with only
finite many
 singular  points $\{q_{j}\}$ in the Gromov-Hausdorff sense.
Furthermore, on $X\backslash \{q_{j}\}$,  $h$ is a  K\"ahler
metric satisfying
 the K\"ahler-Ricci soliton equation.
\end{theorem}
\vskip 5mm

The organization of the paper is as follows: In $\S$2,  we give an
estimate for the  harmonic radius of a solution $g(t)$  of the
normalized K\"ahler-Ricci flow (\ref{1.1}), which plays a central
role in the proof of Theorem 1.1. Then we prove Theorem 1.1 in
$\S$3.

\vskip 10mm

\noindent {\bf Acknowledgement:} The second author   thanks Zhenlei
Zhang and  Professor Xiaochun Rong  for  some helpful discussions.

\section{Main Estimates }

\vskip 4mm

If $(M, g)$ is a complete Riemannian manifold, for a fixed
constant $\Lambda >1$, the {\it harmonic radius} $r_{h}(g)(x)$ at
$x\in M$ is the largest radius of a geodesic ball centered at $x$
on which there are harmonic coordinates $\{h_{i}\}$  such that, if
$g_{ij}=g(\nabla h_{i}, \nabla h_{j})$, then $g_{ij}=\delta_{ij}$
and
$$\Lambda^{-1}\cdot I\leq (g_{ij}) \leq \Lambda \cdot I,$$
$$r_{h}(g)^{1+\alpha}\|g_{ij}\|_{C^{1,\alpha}}\leq \Lambda ,$$ on $B_{g}(x,
r_{h}(g)(x))$ (cf. \cite{An}). In \cite{An}, a lower bound of
$r_{h}(g)$ is obtained by assuming a uniform bound for Ricci
curvature, and a small $L^{n}$-norm for curvature operator.  The
goal of this section is to generalize this result to the solution of
the normalized K\"ahler-Ricci flow (\ref{1.1}).

Let $g(t)$, $t\in [0, +\infty)$, be a solution of the normalized
K\"ahler-Ricci flow (\ref{1.1}) on a compact K\"ahler $n$-manifold
$M$ with $c_{1}(M)>0$ such that $g (0)\in 2\pi c_{1}(M)$.  Recall
that Perelman has shown that there are constants $C$ and  $\kappa$
independent of $t$ such that
\begin{equation}\label{2.1}|R(g(t))|\leq C,\end{equation} and $g(t)$
is $\kappa$-noncollapsed, i.e. for any $r< 1$,  if $|R(g(t))|\leq
r^{-2}$ on a metric ball $B_{g(t)}(x,r)$, then
\begin{equation}\label{2.2}\text{Vol}_{g(t)}(B_{g(t)}(x,r))\geq \kappa r^{2n}, \end{equation} (cf. \cite{ST}).

 \begin{proposition} There are   constants $\varepsilon >0$ and
$\overline{C}>0$ independent of $t$ such that, for any $t\in [1,
+\infty)$,  if, on a metric ball $B_{g(t)}(x,2r)$,
$$\int_{B_{g(t)}(x,2r)}|Rm(g(t))|^{n}dv_{t}\leq \varepsilon,$$ then the
 harmonic radius $r_{h}(g(t))$ satisfies
$$\inf\limits_{B_{g(t)}(x,\frac{r}{2})}r_{h}(g(t)) \geq
\overline{C}r.$$
\end{proposition}

\begin{proof} Actually, we should prove that there are   constants $\varepsilon >0$ and
$\overline{C}>0$ independent of $t$ such that, for any $t\in [1,
+\infty)$,  if, on a metric ball $B_{g(t)}(x,2r)$,
\begin{equation}\label{2.3}\int_{B_{g(t)}(x,2r)}|Rm(g(t))|^{n}dv_{t}\leq \varepsilon,\end{equation} then
\begin{equation}\label{2.4}\frac{r_{h}(g(t))(y) }{d_{g(t)}(y, \partial B_{g(t)}(x,r))}\geq
\overline{C},\end{equation} where $d_{g(t)}(y, \partial
B_{g(t)}(x,r))=\text{dist}_{g(t)}(y, \partial B_{g(t)}(x,r))$. If
it is not true, there is a sequence of times $\{t_{k}\}$, and two
sequences of points $\{x_{k}\}$, $\{y_{k}\}$ such that
$$\frac{r_{h}(g_{k})(y)}{d_{g_{k}}(y,
\partial B_{k})}\geq \frac{r_{h}(g_{k})(y_{k})}{d_{g_{k}}(y_{k},
\partial B_{k})}\longrightarrow
0,$$ for all  $y\in B_{k}$,   when $k\longrightarrow \infty$, where
$g_{k}=g(t_{k})$ and $B_{k}= B_{g_{k}}(x_{k},r)$, but
$$\int_{B_{g_{k}}(x_{k},2r)}|Rm(g_{k})|^{n}dv_{k}\leq \varepsilon.$$

If  $\mu_{k}=r_{h}^{-2}(g_{k})(y_{k})$ and $\tilde{g}_{k}=\mu_{k}
g_{k}$, then $r_{h}(\tilde{g}_{k})(y_{k})=1$,
$$|R(\tilde{g}_{k})|=\mu_{k}^{-1}|R(g_{k})|\leq \frac{
C }{\mu_{k}}\longrightarrow 0, \ \ \rm by \ \ \   (\ref{2.1}),  \ \
\ and $$ $$ d_{\tilde{g}_{k}}^{2}(y_{k},
\partial B_{k})=d_{g_{k}}^{2}(y_{k},
\partial B_{k})\mu_{k}\longrightarrow
 \infty.$$  Furthermore, for any finite  $\rho >0$ and  $z\in
 B_{k}$, if
$d_{\tilde{g}_{k}}(z, y_{k})=dist_{\tilde{g}_{k}}(z, y_{k})< \rho$,
\begin{equation}\label{2.5}r_{h}(\tilde{g}_{k})(z)\geq
\frac{d_{\tilde{g}_{k}}(z,
\partial B_{k})}{d_{\tilde{g}_{k}}(y_{k},
\partial B_{k})}\geq \frac{d_{\tilde{g}_{k}}(z,
\partial B_{k})}{d_{\tilde{g}_{k}}(z,
\partial B_{k})+d_{\tilde{g}_{k}}(z, y_{k})}\geq  \frac{1}{2},\end{equation} $k\gg
1$. Note that, by (\ref{2.2}),
$\text{Vol}_{\tilde{g}_{k}}(B_{\tilde{g}_{k}}(z,\rho))\geq \kappa
\rho^{2n}$ where $z\in B_{k}$, and $\rho < \mu_{k}^{\frac{1}{2}}
\min\{1, C^{-\frac{1}{2}}\} $. Thus, by Lemma 2.1 and Remark 2.4
in \cite{An},  a subsequence of $(B_{k}, \tilde{g}_{k}, y_{k})$
converges to a complete $C^{1,\alpha}$-Riemannian manifold $(N,
g_{\infty}, y_{\infty})$ in the $C^{1,\alpha'}$-sense,
$\alpha'<\alpha$,  i.e. for any $\bar{r}>0$ and $k\gg 1$, there is
a smooth embedding $F_{\bar{r},k}: B_{g_{\infty}}(y_{\infty},
\bar{r}+1) \longrightarrow B_{k}$ such that
$F_{\bar{r},k}^{*}\tilde{g}_{k}$ converges to $g_{\infty}$ in the
$C^{1,\alpha'}$ sense.
 Furthermore,
\begin{equation}\label{2.6} r_{h}(g_{\infty})(y_{\infty})\leq
 \liminf\limits_{k\rightarrow \infty}
 r_{h}(\tilde{g}_{k})(y_{k})=1  \end{equation}  (c.f. \cite{An}), and \begin{equation}\label{2.7}
  \text{Vol}_{g_{\infty}}(B_{g_{\infty}}(y_{\infty},\rho))\geq \kappa
\rho^{2n}\ \ \  \end{equation} for any $\rho >0 $  by (\ref{2.2}).

 \begin{claim} $g_{\infty}$ is a $L^{2,n}$-metric on $N$,
 and a subsequence of  $F_{\bar{r},k}^{*}\tilde{g}_{k}$ converges to $g_{\infty}$ in
 the weak
$L^{2,n}$-topology. Thus,  \begin{equation}\label{2.8}
\int_{N}|Rm(g_{\infty})|^{n}dv_{\infty} \leq \sup\limits_{k}
 \int_{B_{g_{k}}(x_{k},2r)}|Rm(g_{k})|^{n}dv_{k}\leq \varepsilon.
\end{equation}
 \end{claim}

 \begin{proof} Recall that, for a  Riemannian metric $g$, the Ricci
 curvature in harmonic coordinates is given $$g^{ij}\frac{\partial^{2}g_{hl}}{\partial x_{i}\partial x_{j}}
 + \mathcal{Q}(\frac{\partial g_{rs}}{\partial x_{m}})= (Ric(g))_{hl},
 $$ (c.f. \cite{An})  where $\mathcal{Q}$ is a quadratic term.  Note that this
 equation is a uniformly elliptic system of P.D.E. with a  uniform
 $C^{1,\alpha'}$-bound  on the coefficients $g^{ij}$, and a
 $C^{0,\alpha'}$-bound
 on the term $\mathcal{Q}$. If there is a $L^{n}$-bound  on the right
 side, then
 the elliptic regularity theory gives a uniform bound  on
 $\|g\|_{L^{2,n}}$ (c.f. \cite{GT}) . Thus the conclusion follows.
 \end{proof}

   Note that, for   any $\bar{r}>0$,  there is an $r_{0}>0$ such that,
 for $z\in B_{g_{\infty}}(y_{\infty}, \bar{r}+1)$ and  any domain
$\Omega\subset B_{g_{\infty}}(z, r_{0})$,
$\text{Vol}_{g_{\infty}}(\partial \Omega)^{2n}\geq
(1-\delta)c_{2n}\text{Vol}_{g_{\infty}}( \Omega)^{2n-1}$, where
$\delta$ is the constant in Theorem 10.1 in \cite{P1} or Theorem
29.1 of \cite{KL}, and $c_{2n}$ is the Euclidean isoperimetric
constant. Thus, from  the convergence,  for any $z'\in
B_{\tilde{g}_{k}}(y_{k}, \bar{r}+1)$ and any domain $\Omega\subset
B_{\tilde{g}_{k}}(z', r_{0})$,
$\text{Vol}_{\tilde{g}_{k}}(\partial \Omega)^{2n}\geq
(1-\delta)c_{2n}\text{Vol}_{\tilde{g}_{k}}( \Omega)^{2n-1}$, when
$k\gg 1$.

\begin{lemma} $g_{\infty}$ is a complete  Ricci flat metric on $N$.
\end{lemma}

\begin{proof}
 Note that  $\tilde{g}_{k}(t)=\mu_{k}g(t_{k}+\mu_{k}^{-1}t)$, $t\in
[0, \infty)$, are solutions of the normalized K\"ahler-Ricci flow
 $$ \partial_{t}\tilde{g}_{k}(t)=-Ric(\tilde{g}_{k}(t))+\frac{1}{\mu_{k}}\tilde{g}_{k}(t)$$ with initial metrics  $\tilde{g}_{k}$, which satisfy
\begin{equation}\label{2.9}|R(\tilde{g}_{k}(t))|=\mu_{k}^{-1}|R(g(t_{k}+\mu_{k}^{-1}t))|\leq
\frac{ C}{\mu_{k}}\longrightarrow 0,  \ \ \rm by \ \ \
(\ref{2.1}),\end{equation} when $k\longrightarrow \infty$.  If
$\bar{t}=\mu_{k}(1-e^{-\frac{t}{\mu_{k}}})$ and
 $\bar{g}_{k}(\bar{t})=e^{-\frac{t}{\mu_{k}}}\tilde{g}_{k}(t)$, then $\bar{g}_{k}(\bar{t})$, $\bar{t}\in [0, \mu_{k})$ is a
 solution of the  K\"ahler-Ricci
flow
$$\partial_{\bar{t}}\bar{g}_{k}(\bar{t})=-Ric(\bar{g}_{k}(\bar{t}))$$
on $M$ with initial metric $\tilde{g}_{k}$.   By Theorem 10.1 in
\cite{P1} or Theorem 29.1 in \cite{KL}, there is an $\epsilon >0$
such that
\begin{equation}\label{2.10}|Rm(\bar{g}_{k}(\bar{t}))|(z)\leq
\frac{1}{\bar{t}}+(\epsilon r_{0})^{-2},\end{equation}
$0<\bar{t}<(\epsilon r_{0})^{2}<1$,  for all
 $z\in
B_{\tilde{g}_{k}}(y_{k}, \bar{r}+1)$. By (\ref{2.9}),
\begin{equation}\label{2.11}|R(\bar{g}_{k}(\bar{t}))|=\frac{1}{1-\frac{\bar{t}}{\mu_{k}}}|R(\tilde{g}_{k}(t))|\leq
\frac{2 C}{\mu_{k}}\longrightarrow 0, \end{equation} $\bar{t}\in [0,
(\epsilon r_{0})^{2}]$, when $k\longrightarrow \infty $.

 By (\ref{2.2}), it is easy to see that
  $\bar{g}_{k}(\bar{t})$ is  $\kappa$-noncollapsed, and, for any $z\in
B_{\tilde{g}_{k}}(y_{k}, \bar{r}+1)$, the injectivity radius
$\text{inj}_{\bar{g}_{k}((\epsilon r_{0})^{2} )}(z)\geq \iota$ for
a positive constant $\iota$ independent of $k$ by (\ref{2.10}).
From the compactness theorem for Ricci flow (c.f. Appendix E in
\cite{KL} or \cite{H1}), by passing to a subsequence,
$(B_{\tilde{g}_{k}}(y_{k}, \bar{r}-1), \bar{g}_{k}(\bar{t}),
y_{k})$, $\bar{t}\in (0, T]$, $C^{\infty}$-converges to
$(B_{\infty}, g_{\infty}(\bar{t}), y_{\infty})$, $\bar{t}\in (0,
T]$, where $g_{\infty}(\bar{t})$ is a solution of Ricci flow on
$B_{\infty}$, and $T< (\epsilon r_{0})^{2}$. By (\ref{2.11}),
$|R(g_{\infty}(\bar{t}))|\equiv 0$. This implies that
$|Ric(g_{\infty}(\bar{t}))|\equiv 0$, and
$g_{\infty}(\bar{t})\equiv g_{\infty}(T)$  is a Ricci-flat metric
on $B_{\infty}$. By Theorem 36.2 in \cite{KL},
$$d_{GH}((B_{\infty}, g_{\infty}(T)), (B_{g_{\infty}}(y_{\infty}, \bar{r}-1),
g_{\infty}))= \lim\limits_{\bar{t}\rightarrow 0}d_{GH}((B_{\infty},
g_{\infty}(\bar{t})), (B_{g_{\infty}}(y_{\infty}, \bar{r}-1),
g_{\infty}))=0,$$ where $d_{GH}$ denotes the Gromov-Hausdorff
distance.  By letting $ \bar{r}\rightarrow\infty$ and taking a
diagonalized sequence, we obtain that $g_{\infty}$ is a Ricci flat
metric on $N$.
\end{proof}

 Lemma  2.3 and (\ref{2.7})  imply
that there is a global bound for the Sobolev constant on $(N,
g_{\infty})$ (c.f. \cite{Cr} or \cite{An1}). Let $\varepsilon$ be
the corresponding constant  in Lemma 2.1 of \cite{An1}, which
depends on the Sobolev constant on $(N, g_{\infty})$.  By Claim 2.2,
 Lemma 2.3 and  Lemma 2.1 in \cite{An1}, $$\sup_{B_{g_{\infty}}(y_{\infty},
\frac{s}{2})}|Rm(g_{\infty})|\leq \frac{C'}{s^{2}} \longrightarrow
0, \ \ \ \ {\rm as} \ \   \ s\rightarrow\infty,$$ and, thus,
$g_{\infty}$ is a  flat  metric. By  (\ref{2.7}),  $(N,
g_{\infty})$ is the standard Euclidean space $\mathbb{R}^{2n}$
(c.f. \cite{An1}). It contradicts to (\ref{2.6}), since  the
harmonic radius of $\mathbb{R}^{2n}$  is infinite. We obtain the
conclusion.
\end{proof}

\vskip 4mm

\section{Proof of Theorem 1.1}

\vskip 4mm

Let $g(t)$, $t\in [0, +\infty)$, be a solution of the normalized
K\"ahler-Ricci flow (\ref{1.1}) on a compact K\"ahler $n$-manifold
$M$ with $c_{1}(M)>0$ such that $g (0)\in 2\pi c_{1}(M)$, and
\begin{equation}\label{3.1} \int_{M}|Rm(g(t))|^{n}dv_{t}\leq C,\end{equation}
where $C$ is a constant independent of $t$. Assume that
$t_{k}\longrightarrow \infty$ is a sequence of times.

\vskip 4mm

\begin{lemma} If there is a constant $\mathcal{V}>0$ independent of
$t$  such that
\begin{equation}\label{3.2}\text{Vol}_{g(t)}(B_{g(t)}(x, r))\leq
\mathcal{V} r^{2n},
\end{equation} for any $r\leq 1$ and  $x\in M$, then, by passing to a subsequence,
$(M, g(t_{k}))$ converges to an orbifold $(X, h)$ with only
finitely many singular  points $\{q_{j}\}$ in the Gromov-Hausdorff
topology,  where $h$ is a $C^{0}$-orbifold metric, and is $C^{1,
\alpha}$ off the
 singular  points. Furthermore, for any compact subset $K\subset
X\backslash \{q_{j}\}$, there are smooth embeddings $F_{K,k}:
K\longrightarrow M$ such that $F_{K,k}^{*}g(t_{k})$ converges to
$h$ in the $C^{1,\alpha'}$ (resp. weak $L^{2,n}$) topology,
$\alpha'<\alpha$.
\end{lemma}

\begin{proof} By (\ref{2.2}) and (\ref{3.1}), all  hypothesis  of
Theorem 1.1 in \cite{An2} are satisfied except (1.5) in \cite{An2},
i.e. small curvature estimate, which is used in the arguments in
Section 2.1 of \cite{An2} for obtaining the local
$C^{1,\alpha}$-convergence, and weak $L^{2,p}$-convergence. However,
we have Proposition 2.1, which gives a lower bound for harmonic
radius, when the $n$-norm of curvature is small enough. This is
enough to obtain the local $C^{1,\alpha'}$-convergence (c.f.
\cite{An}), and weak $L^{2,n}$-convergence by the proof of Claim
2.2.  Thus, by (\ref{3.2}),  the arguments in Section 2.1 of
\cite{An2} and the proof of Theorem 2.6
 in  \cite{An}, we obtain that, by passing to a subsequence,
$(M, g(t_{k}))$ converges to a metric space  $(X, h)$ in the
Gromov-Hausdorff topology.  Here $X$ is a multi-fold  with a finite
set  $S=\{q_{j}\}$  of
 singular  points in the sense of \cite{TV}, i.e. $X$ is a smooth
manifold off a finite set of singular points $S$, and a neighborhood
of each singular  point $q_{j}$  is a finite union of  cones  on
  spherical space forms  $C(S^{2n-1}/\Gamma)$, $\Gamma\subset SO(2n)$, where the vertex of each
  cone is identified with the point $q_{j}$. The
metric $h$  is a $C^{1,\alpha}\cap L^{2,n}$-Riemannian metric on
$X\backslash S$, and, in a local uniformization $B^{2n}\backslash
\{0\}$, $h$ extends to a $C^{0}$-metric on the ball $B^{2n}$.
Furthermore, for any compact subset $K\subset X\backslash S$, there
are smooth embeddings $F_{K,k}: K\longrightarrow M$ such that
$F_{K,k}^{*}g_{k}$ converges to $h$ in the $C^{1,\alpha'}$ and weak
$L^{2,n}$ topologies, $\alpha'<\alpha$, where $g_{k}=g(t_{k})$. Thus
the only thing we are supposed to prove is that there is exactly one
cone at each singular point.

From the arguments in Section 2.1 of \cite{An2}, for any singular
point $x\in S$, there is a sequence of points $x_{k} \in ( M, g_{k})
$ such that $x_{k}\longrightarrow x$, when $k\rightarrow\infty$,
and, for any $r>0$,
\begin{equation}\label{3.3}\liminf\limits_{k\rightarrow\infty}\int_{B_{g_{k}}(x_{k},
r)}|Rm(g_{k})|^{n}dv_{k}\geq \varepsilon,\end{equation} where
$\varepsilon$ is the constant in Proposition 2.1. If there is a
$x\in S$ with more than one  cones attaching  it,  we imitate the
proof in \cite{TV} to obtain a contradiction.  Choose a radius
$\bar{r}>0$ small enough, and a sequence of points $x_{k}\in ( M,
g_{k}) $ such that the harmonic radius  satisfies
\begin{equation}\label{3.4}\inf\limits_{B_{g_{k}}(x_{k},
\bar{r})}i_{h}(g_{k})=i_{h}(g_{k})(x_{k}) \longrightarrow 0, \ \ \ \
{\rm as} \ \ \ \ k\rightarrow \infty,
\end{equation} and \begin{equation}\label{3.5}\int_{B_{h}(x,
\bar{r})}|Rm(h)|^{n}dv_{h}\leq \frac{\varepsilon}{2},
\end{equation} and $B_{h}(x,
\bar{r})\backslash\{x\}$ has several components.  We choose $r_{k}<
\bar{r}$ such that
\begin{equation}\label{3.6} \int_{D_{g_{k}}(r_{k},
\bar{r})}|Rm(g_{k})|^{n}dv_{k}=\varepsilon, \end{equation} where
$D_{g_{k}}(r_{k}, \bar{r})=B_{g_{k}}(x_{k}, \bar{r}) \backslash
B_{g_{k}}(x_{k}, r_{k})$. Note that the annulus $D_{g_{k}}(r_{k},
\bar{r})$ has several components.

  By (\ref{3.4}) and  (\ref{3.5}), it is easy to see that   $x_{k}\longrightarrow x$ and
$r_{k}\longrightarrow 0$ when $k\rightarrow\infty $. From the same
arguments as above, a subsequence of $(M, r_{k}^{-2}g_{k}, x_{k})$
converges to a complete multi-fold $(N_{1}, g_{1,\infty},
x_{1,\infty})$  with a finite set $S_{1}$ of  singular  points in
the pointed Gromov-Hausdorff sense, where $ g_{1,\infty}$ is a
$C^{1,\alpha}\cap L^{2,n}$-Riemannian metric on $N_{1}\backslash
S_{1}$. And, for any compact subset $K\subset N_{1}\backslash
S_{1}$, there are smooth embeddings $F_{1,K,k}: K\longrightarrow M$
such that $F_{1,K,k}^{*}r_{k}^{-2}g_{k}$ converges to $g_{1,\infty}$
in the $C^{1,\alpha'}$ and weak $L^{2,n}$ topologies,
$\alpha'<\alpha$. By (\ref{3.6}), there is no singular points in
$N_{1}\backslash B_{g_{1,\infty}}(x_{1,\infty}, 1)$.  By considering
the rescaled Ricci flow $r_{k}^{-2}g(t_{k}+r_{k}^{2}t)$, $t\in [0,
\infty)$,  the arguments in the proof of Lemma  2.3 imply   that $
g_{1,\infty}$ is a Ricci-flat  metric on $N_{1}\backslash S_{1}$.

If, for any $x_{1}\in S_{1}$, there is only one cone  at $x_{1}$,
i.e. $N_{1}$ is an orbifold, then $g_{1,\infty}$ is a Ricci flat
orbifold  metric on $N_{1}$ by (\ref{3.1}) and \cite{An1}, i.e.  in
a local uniformization $B^{2n}\backslash \{x_{1}\}$, $g_{1,\infty}$
extends to a  Ricci flat metric on the ball $B^{2n}$. From the
construction above, $N_{1}$  has several ends, i.e. $N_{1}\backslash
B_{g_{1,\infty}}(x_{1,\infty}, 1)$ is not connected. By the
convergence and the  Perelman's estimate (\ref{2.2}),  there is a
$\kappa
>0$ such that, for any $y_{1}\in N_{1}$ and $\rho >0$,
$$Vol_{g_{1,\infty}}(B_{g_{1,\infty}}(y_{1},\rho))\geq \kappa \rho^{2n}.$$
 This implies that   there is a global bound for the Sobolev constant on
$N_{1}\backslash B_{g_{1,\infty}}(x_{1,\infty}, 1)$ (c.f. \cite{Cr}
or \cite{An1}).  Thus, by Lemma 2.1 of \cite{An1}, the small
curvature estimate (1.5) in \cite{An2} is satisfied. Theorem 1.2 and
Remark 2.8 of \cite{An2} show  that $(N_{1}, g_{1,\infty},
x_{1,\infty})$ is an Asymptotically Locally Euclidean multi-fold
with several ends. Each end $E_{j}$ of $N_{1}$ is diffeomorphic to
$(\mathbb{R}^{2n}\backslash B^{2n})/\Gamma_{j}$, where $\Gamma_{j}$
is a finite subgroup of $SO(2n)$ acting on $S^{2n-1}$ freely.
However, since the number of   ends of $N_{1}$ is larger than one,
the splitting theorem for orbifolds  (c.f. \cite{B2}) implies that
$(N_{1}, g_{1,\infty})$ is isometric to $Y\times \mathbb{R}^{1}$
where $Y$ is a complete Ricci-flat orbifold.  It is a contradiction.
Thus there is a $x_{1}\in S_{1}$ such that there are more than one
cones at $x_{1}$.

Now we do the same process as above for $N_{1}$ and $x_{1}$, and
obtain a multi-fold $N_{2}$ with finite singular set $S_{2}$. If
there is only one cone at each point of $S_{2}$, then we stop, and
obtain a contradiction as above. Otherwise we repeat the procedure.
This process must terminate in finite steps, since each singularity
takes at least $\varepsilon$ of curvature by the construction (c.f.
 \cite{TV}).  Finally, we shall obtain a contradiction. Thus  $X$ is an orbifold.
\end{proof}

\begin{lemma}  There is a constant $\mathcal{V}>0$ independent of
$t$ such that \begin{equation}\label{3.7}Vol_{g(t)}(B_{g(t)}(x,
r))\leq \mathcal{V} r^{2n}, \ \ \
\end{equation} for any $r\leq 1$ and $x\in M$.
\end{lemma}

\begin{proof} If it is  not true, there exists a sequence of times $t_{i}$, and a sequence
of balls $B_{g_{i}}(x_{i},r_{i})$ such that
$$\frac{Vol_{g_{i}}(B_{g_{i}}(x_{i},r_{i}))}{r_{i}^{2n}}\geq
\mathcal{V}_{i} \longrightarrow \infty, $$ as $i\rightarrow\infty$,
where $g_{i}=g(t_{i})$. Let $s_{i}$ be the smallest radius such
that, for some $y_{i}\in (M, g_{i})$,
$$Vol_{g_{i}}(B_{g_{i}}(y_{i},s_{i}))\geq 2\omega_{2n}s_{i}^{2n},$$
where $\omega_{2n}$ is the volume of 1-ball in the Euclidean space
$\mathbb{R}^{2n}$. (Note that, for any $x\in M$,
$\text{Vol}_{g_{i}}(B_{g_{i}}(x,r))\sim \omega_{2n}r^{2n}$ when
$r\rightarrow 0$.)  From the arguments in the proof of Theorem 1.1
in  \cite{An2}, we can assume that $s_{i}\longrightarrow 0$ as
$i\rightarrow \infty$.

If  $\tilde{g}_{i}=s_{i}^{-2}g_{i}$, then, for any $r< 1$ and $x\in
M$,
\begin{equation}\label{3.8}Vol_{\tilde{g}_{i}}(B_{\tilde{g}_{i}}(x,r))\leq
2\omega_{2n}r^{2n}, \ \ \ {\rm and} \ \
Vol_{\tilde{g}_{i}}(B_{\tilde{g}_{i}}( y_{i},1))=
2\omega_{2n}\end{equation} by the choice of $s_{i}$. By the
arguments in Section 2.1 in \cite{An2} and the proof of Lemma 3.1, a
subsequence of $(M, \tilde{g}_{i}, y_{i})$ converges to a complete
orbifold $(N, \tilde{g}_{\infty}, y_{\infty})$ with only finite many
 singular  points $\{q_{j}\}$, where $\tilde{g}_{\infty}$ is a
$C^{0}$-orbifold metric, and is $C^{1, \alpha}$ off the singular
points. Furthermore, for any compact subset $K\subset N\backslash
\cup_{j}\{q_{j}\}$, there are smooth embeddings $F_{K,i}:
K\longrightarrow M$ such that $F_{K,i}^{*}\tilde{g}_{i}$ converges
to $\tilde{g}_{\infty}$ in the $C^{1,\alpha'}$ and weak $L^{2,n}$
topologies, $\alpha'< \alpha$. By considering the re-scaled
K\"ahler-Ricci flow
$\tilde{g}_{k}(t)=s_{i}^{-2}g(t_{i}+s_{i}^{2}t)$, $t\in [0,
\infty)$, and the proof of Lemma  2.3, we obtain that
$\tilde{g}_{\infty}$ is a Ricci flat metric on $N\backslash
\cup_{j}\{q_{j}\}$. From the weak $L^{2,n}$-convergence,
$$\int_{N\backslash
\cup_{j}\{q_{j}\}}|Rm(\tilde{g}_{\infty})|^{n}dv_{\infty}\leq
C<\infty.$$ Thus $\tilde{g}_{\infty}$ can be extended to a Ricci
flat orbifold  metric on $N$ (c.f. \cite{An1}). Note that volume
comparison theorem valids also for orbifolds by \cite{B}. Thus,
for any $r>0$, we obtain that
$$\text{Vol}_{\tilde{g}_{\infty}}(B_{\tilde{g}_{\infty}}(x,r))\leq
\omega_{2n}r^{2n}. $$ By (\ref{3.8}) and the
$C^{1,\alpha'}$-convergence, we obtain
$Vol_{\tilde{g}_{\infty}}(B_{\tilde{g}_{\infty}}( y_{\infty},1))=
2\omega_{2n}$, which  is a contradiction.
\end{proof}

\begin{proof}[Proof of Theorem 1.1] By Lemmas 3.1 and 3.2,   a
subsequence of  $(M, g(t_{k}))$ converges to an orbifold $(X, h)$
with  only  finite many  singular  points $\{q_{j}\}$ in the
Gromov-Hausdorff topology,  where $h$ is a $C^{0}$-orbifold metric,
and is $C^{1, \alpha}$ off the singular  points. Furthermore, for
any $r>0$, there are smooth embeddings $F_{r,k}:  X\backslash
\cup_{j}B_{h}(q_{j}, r) \longrightarrow M$ such that
$F_{r,k}^{*}g(t_{k})$ converges to $h$ in the $C^{1,\alpha'}$ and
weak $L^{2,n}$ topologies, $\alpha'<\alpha$. Actually, $h$ is a
K\"ahler metric (c.f. \cite{Ru}). By the same arguments as in the
proof of Proposition 2.1, for $k\gg 1$, there is a $r_{0}$
independent of $k$ such that, for any domain $\Omega\subset
B_{g_{k}}(z, r_{0})$, $z\in X\backslash \cup_{j}B_{h}(q_{j}, 2r) $,
$Vol_{g_{k}}(\partial \Omega)^{2n}\geq (1-\delta)c_{2n}Vol_{g_{k}}(
\Omega)^{2n-1}$, where $g_{k}=g(t_{k})$,  $\delta$ is the constant
in Theorem 29.1 of \cite{KL}, and $c_{2n}$ is the Euclidean
isoperimetric constant.

If $\tilde{t}=1-e^{-t}$ and
 $\tilde{g}_{k}(\tilde{t})=e^{-t}g_{k}(t)$, where $g_{k}(t)=g(t_{k}+t)$,
  then $\tilde{g}_{k}(\tilde{t})$, $\tilde{t}\in [0, 1)$ is a
 solution of the  K\"ahler-Ricci flow
$$\partial_{\tilde{t}}\tilde{g}_{k}(\tilde{t})=-Ric(\tilde{g}_{k}(\tilde{t}))$$
on $M$ with initial metric $g_{k}$.   By  Theorem 10.1 in \cite{P1}
or  Theorem 29.1 in \cite{KL}, there is a $\epsilon >0$ such that
$$|Rm(\tilde{g}_{k}(\tilde{t}))|(z)\leq \frac{1}{\tilde{t}}+(\epsilon r_{0})^{-2},$$ for $0<\tilde{t}<(\epsilon
r_{0})^{2}<1$, and all $z\in  X\backslash \cup_{j}B_{h}(q_{j}, 2r)$.
Since $\tilde{g}_{k}(\tilde{t})$ differs  to $g_{k}(t)$ only by
rescalings in space and time, i.e.
$g_{k}(t)=\frac{1}{1-\tilde{t}}\tilde{g}_{k}(\tilde{t})$,  where
  $t=\log
(\frac{1}{1-\tilde{t}})$, we obtain
\begin{equation}\label{3.9}|Rm(g_{k}(t))|(z)\leq
\frac{1}{e^{-t}(1-e^{-t})}+e^{t}(\epsilon r_{0})^{-2},\end{equation}
for all $z\in X\backslash \cup_{j}B_{h}(q_{j}, 2r)$ and $0<t<-\log
(1-(\epsilon r_{0})^{2})$. By (\ref{2.2}) and (\ref{3.9}),
  $g_{k}(t)$ is  $\kappa$-noncollapsed, and, for any $z\in
X\backslash \cup_{j}B_{h}(q_{j}, 2r)$, the injectivity radius
$\text{inj}_{\bar{g}_{k}(t_{0} )}(z)\geq \iota$, where
$t_{0}=-\frac{1}{2}\log (1-(\epsilon r_{0})^{2})$, for a constant
$\iota$ independent of $k$. By the compactness theorem for Ricci
flow (c.f. Appendix E in \cite{KL} or \cite{H1}), by passing to a
subsequence, $(X\backslash \cup_{j}B_{h}(q_{j}, 4r), g_{k}(t))$,
$t\in (0, T]$, converges to $(X\backslash \cup_{j}B_{h}(q_{j}, 4r),
h(t))$, $t\in (0, T]$, where $h(t)$ is a solution of Ricci flow on
$X\backslash \cup_{j}B_{h}(q_{j}, 4r)$, and   $T< t_{0}$.
 By Theorem 36.2 in \cite{KL},
\begin{equation}\label{3.10} \lim\limits_{t\rightarrow 0}d_{GH}((X\backslash
\cup_{j}B_{h}(q_{j}, 4r), h), (X\backslash \cup_{j}B_{h}(q_{j}, 4r),
h(t)))=0.\end{equation}

From the arguments in  the proof of Theorem 12 in \cite{ST}, $h(t)$,
$t\in (0,T]$, satisfies the K\"ahler-Ricci soliton equation, i.e.
there are smooth functions $f(t)$, $t\in (0,T]$, such that
$$Ric(h(t))-h(t)=\nabla \overline{\nabla}f(t),$$ and
$$\nabla \nabla f(t)=\overline{\nabla} \overline{\nabla} f(t)=0.
$$ Note that there is a fixed vector field $V$ on $X\backslash \cup_{j}B_{h}(q_{j}, 4r)$  such that $$h(t)=\phi^{-1}(t)^{*}h(t_{1}),
 \ \ \ \ {\rm and } \ \ \ f(t)=\phi^{-1}(t)^{*}f(t_{1}),$$ where  $\{\phi(t)\}$   is the 1-parameter
 group of diffeomorphisms generated by $-V$, and $\phi(t_{1})=id$ (c.f. Appendix C of
 \cite{KL}). By letting $0< t_{1}\ll 1$ such that $\phi^{-1}(0)(X\backslash \cup_{j}B_{h}(q_{j}, 8r)
 )\subset X\backslash \cup_{j}B_{h}(q_{j}, 4r)$, we obtain that
 $\phi^{-1}(0) ^{*}h(t_{1})$ is a metric satisfying  the K\"ahler-Ricci soliton
 equation on $X\backslash \cup_{j}B_{h}(q_{j}, 8r)$, and $$\lim\limits_{t\rightarrow 0}d_{GH}((X\backslash
\cup_{j}B_{h}(q_{j}, 8r), \phi^{-1}(0) ^{*}h(t_{1})), (X\backslash
\cup_{j}B_{h}(q_{j}, 8r), h(t)))=0.$$ Thus, by (\ref{3.10}),
$(X\backslash \cup_{j}B_{h}(q_{j}, 8r), h)$ is isometric to
$(X\backslash \cup_{j}B_{h}(q_{j}, 8r), \phi^{-1}(0) ^{*}h(t_{1}))$.
By letting $r\rightarrow 0$ and taking  a diagonalized sequence, we
obtain that $h$ satisfies   the K\"ahler-Ricci soliton
 equation on $X\backslash \cup_{j}\{q_{j}\}$.

\end{proof}

\end{document}